\documentclass{article}

\usepackage{amsmath,amsfonts,amssymb,url,fullpage,color,chngcntr,MnSymbol}
\usepackage[utf8]{inputenc}

\usepackage{theorem}
\usepackage{graphicx}

{\theorembodyfont{\slshape}
\newtheorem{proposition}{Proposition}[section]
\newtheorem{lemma}[proposition]{Lemma}
\newtheorem{corollary}[proposition]{Corollary}
\newtheorem{theorem}[proposition]{Theorem}}

{\theorembodyfont{\upshape}
\newtheorem{remark}[proposition]{Remark}
}

\newcommand\SL{{\mathop{\mathrm{SL}}}}
\newcommand\ricf{{\mathop{\mathrm{RiCF}}}}

\newcommand\End{{\mathop{\mathrm{End}}}}
\newcommand\LCM{{\mathop{\mathrm{LCM}}}}
\newcommand\Gal{{\mathop{\mathrm{Gal}}}}

\renewcommand\Im{{\mathop{\mathrm{Im}}}}

\newcommand\gal{{\mathop{\mathrm{Gal}}}}
\newcommand\Cl{{\mathop{\mathrm{Cl}}}}

\newcommand\C{{\mathbb C}}

\newcommand\Q{{\mathbb Q}}

\newcommand\Z{{\mathbb Z}}

\newcommand\OO{{\mathcal O}}

\newcommand\eps\varepsilon

\newcommand\tilN{{\widetilde N}}

\newcommand\qed{\hfill$\square$}

\title{CM-Points on Straight Lines}
\author{Bill Allombert, Yuri Bilu, Amalia Pizarro-Madariaga}
\date\today

1

\makeatletter

\renewcommand*\l@section[2]{%
  \ifnum \c@tocdepth >\z@
    \addpenalty\@secpenalty
    \addvspace{0.3em \@plus\p@}%
    \setlength\@tempdima{1.5em}%
    \begingroup
      \parindent \z@ \rightskip \@pnumwidth
      \parfillskip -\@pnumwidth
      \leavevmode \bfseries
      \advance\leftskip\@tempdima
      \hskip -\leftskip
      #1\nobreak\hfil \nobreak\hb@xt@\@pnumwidth{\hss #2}\par
    \endgroup
  \fi}
  
 \makeatother

\counterwithin{table}{section}

\begin{document}

\hfuzz=5pt

\maketitle

\begin{abstract}
We prove that, with ``obvious'' exceptions, a CM-point $(j(\tau_1),j(\tau_2))$ cannot belong to a straight line in~$\C^2$ defined over~$\Q$. This generalizes a result of Kühne, who proved this for the line ${x_1+x_2=1}$. 
\end{abstract}

{\footnotesize

\tableofcontents

}

\section{Introduction}

In this article~$\tau$ with or without indices denotes a  quadratic\footnote{``Quadratic'' here and below mean``of degree~$2$ over~$\Q$~''.} complex number with ${\Im\tau>0}$ and~$j$ denotes the $j$-invariant.

In 1998 André~\cite{An98} proved that a non-special irreducible plane curve in $\C^2$ may have only finitely many CM-points. Here a \textsl{plane curve} is a curve defined by an irreducible polynomial equation ${F(x_1,x_2)=0}$, where~$F$ is a polynomial with complex coefficients. A \textsl{CM-point} in $\C^2$ is a point of the form $(j(\tau_1),j(\tau_2))$ with  quadratic $\tau_1,\tau_2$. \textsl{Special curves} are the curves of the following types:

\begin{itemize}
\item
``vertical lines'' ${x_1=j(\tau_1)}$;
\item
``horizontal lines'' ${x_2=j(\tau_2)}$;

\item
\textsl{modular curves} $Y_0(N)$, realized as the plane curves ${\Phi_N(x_1,x_2)=0}$, where $\Phi_N$ is the modular polynomial of level~$N$.
\end{itemize}

Clearly, each special curve contains infinitely many CM-points, and André proved that special curves are characterized by this property.

André's result was the first non-trivial contribution to the celebrated André-Oort conjecture on the special subvarieties of Shimura varieties; see~\cite{Pi11} and the references therein. 

Independently of André the same result was also obtained by Edixhoven~\cite{Ed98}, but Edixhoven had to assume the Generalized Riemann Hypothesis for certain $L$-series to be true. 

Further proof followed; we mention specially the remarkable argument of Pi\-la~\cite{Pi09}. It is based on an idea of Pila and Zannier~\cite{PZ08} and readily extends to higher dimensions~\cite{Pi11}.

The arguments mentioned above were non-effective, because they used the Siegel-Brauer lower bound for the class number. Breuer~\cite{Br01} gave an effective proof, but it depended on GRH.

Recently Kühne~\cite{Ku12,Ku13} and, independently, Bilu, Masser, and Zannier~\cite{BMZ13} found unconditional effective proofs of André's theorem. Besides giving general results, both articles~\cite{Ku13} and~\cite{BMZ13} treat also some particular curves, showing they have no CM-points at all. For instance, Kühne \cite[Theorem~5]{Ku13} proves the following.

\begin{theorem}
\label{thx+y}
The straight line ${x_1+x_2=1}$ has no CM-points. 
\end{theorem}

(The same result was also independently obtained in an earlier version of \cite{BMZ13}, but did not appear in the final version.)

A similar result for the curve ${x_1x_2=1}$ was obtained in \cite{BMZ13}. 

One can ask about CM-points on general straight lines defined over~$\Q$; that is, defined by an equation \begin{equation}
\label{estraightline}
A_1x_1+A_2x_2+B=0,
\end{equation}
where ${A_1,A_2,B\in \Q}$. One has to exclude from consideration the \textsl{special straight lines}: ${x_1=j(\tau_1)}$, ${x_2=j(\tau_2)}$ and ${x_1=x_2}$, the latter being nothing else than the modular curve $Y_0(1)$ (the modular polynomial $\Phi_1$ is ${x_1-x_2}$).  According to the theorem of André, these are the only straight lines containing inifinitely many CM-points.

In the present paper we obtain a rather vast generalization of Theorem~\ref{thx+y}. 

\begin{theorem}
\label{thmain}
Let $(j(\tau_1),j(\tau_2))$ be a CM-point belonging to a non-special straight line defined over~$\Q$. Then we have one of the following options. Either
\begin{equation}
\label{edeg1}
j(\tau_1),j(\tau_2)\in \Q,
\end{equation}
or 
\begin{equation}
\label{edeg2}
j(\tau_1)\ne j(\tau_2), \quad \Q(j(\tau_1))=\Q(j(\tau_2)), \quad [\Q(j(\tau_1)):\Q]=[\Q(j(\tau_2)):\Q]=2.
\end{equation}
\end{theorem}

\begin{remark}

\begin{enumerate}

\item
Recall that ${[\Q(j(\tau):\Q]=h(\OO_\tau)}$, the class number  of the ``complex multiplication order'' ${\OO_\tau=\End\langle\tau,1\rangle}$, where ${\langle\tau,1\rangle}$ is the lattice generated by~$\tau$ and~$1$. All orders of class number~$1$ and~$2$ are well-known, which means that points satisfying~\eqref{edeg1} or~\eqref{edeg2} can be easily listed. In fact, there are 169 CM-points satisfying~\eqref{edeg1} and, up to $\Q$-conjugacy, 217 CM-points satisfying~\eqref{edeg2}; see Remark~\ref{rall} for the details. 

\item
Our result is best possible because any point satisfying~\eqref{edeg1} or~\eqref{edeg2} does belong to a non-special straight line defined over~$\Q$.

\item
Kühne remarks on page~5 of his article~\cite{Ku13} that his Theorem~4 allows one, in principle, to list all possible CM-points belonging to non-special straight lines over~$\Q$, but the implied calculation does not seem to be feasible. 

\item
Bajolet~\cite{Ba14} produced a software package for finding all CM-points on a given straight line. He illustrated its efficiency by proving that no  straight line~\eqref{estraightline} with non-zero ${A_1,A_2,B\in \Z}$ satisfying ${|A_1|,|A_2|,|B|\le 10}$ passes through a CM-point.  This work is now formally obsolete because of our Theorem~\ref{thmain}, but a similar method can be used in more general situations, where our theorem no longer applies. 

\item
CM-points $(x_1,x_2)$ satisfying ${x_1x_2\in \Q^\times}$ are completely classified in~\cite{BLP14}; this generalizes the above-mentioned result from~\cite{BMZ13} about the curve ${x_1x_2=1}$.

\end{enumerate}

\end{remark}

In Sections~\ref{scfi} and~\ref{sricf} we recall basic facts about imaginary quadratic orders, class groups, ring class fields and complex multiplications. 

In Section~\ref{stwo} we investigate the 
field equality ${\Q(j(\tau_1))=\Q(j(\tau_2))}$. In particular, in Corollary~\ref{cojopa} we determine all cases of such equality when ${\Q(\tau_1)\ne \Q(\tau_2)}$. This might be of independent interest.

After all these preparations, we prove Theorem~\ref{thmain} in Section~\ref{sproof}.

\paragraph{Acknowledgments.}
We thank Lars Kühne whose marvelous article~\cite{Ku13} was our principal source of inspiration. We also thank Karim Belabas, Henri Cohen, Andreas Enge  and Jürg Kramer for useful conversations, and the referee for the encouraging report and many helpful comments. 

Yuri Bilu was supported by  the \textsl{Agence National de la Recherche} project ``Hamot'' (ANR 2010 BLAN-0115-01).  Amalia Pizarro-Madariaga was supported by the ALGANT scholarship program.

Our calculations were performed using the \textsf{PARI/GP} package~\cite{PARI}.

\section{Imaginary Quadratic Orders}
\label{scfi}

In this section we recall basic  facts about imaginary quadratic fields and their orders, and recall a famous result of Weinberger about class groups annihilated by~$2$.

\subsection{Class Groups}
\label{ssclgr}

Let~$K$ be an imaginary quadratic field and~$\OO$ an  order in~$K$ of discriminant ${\Delta=Df^2}$, where~$D$ is the discriminant of the field~$K$ (often called the \textsl{fundamental discriminant}) and~$f$ is the conductor of~$\OO$, defined from ${\OO=\Z+f\OO_K}$. We denote by $\Cl(\OO)$ the class group of~$\OO$ (the group of invertible fractional ideals modulo the invertible principal fractional ideals). As usual, we set ${h(\OO)=|\Cl(\OO)|}$. Since~$\OO$ is uniquely determined by its discriminant~$\Delta$, we may also write  $\Cl(\Delta)$, $h(\Delta)$  etc. In particular, ${\Cl(D)=\Cl(\OO_K)}$ is the  class group of the field~$K$, and $h(D)$ is the class number of~$K$. 

There is a canonical exact sequence
\begin{equation}
\label{eexcl}
1\to \Cl_0(\Delta)\to \Cl(\Delta)\to\Cl(D)\to 1,
\end{equation}
where the kernel $\Cl_0(\Delta)$ will be described below. This implies, in particular, that ${h(D)\mid h(\Delta)}$.

The structure of the group $\Cl_0(\Delta)$ is described, for instance, in \cite{Co89}, Section~7.D and Exercise~7.30. We briefly reproduce this description here. We will assume, with a slight abuse of notation, that ${\Z/f\Z}$ is a subring of ${\OO_K/f\OO_K}$. Then have another canonical exact sequence
\begin{equation}
\label{ecl0structure}
1\to (\Z/f\Z)^\times(\OO_K^\times)_f\hookrightarrow (\OO_K/f\OO_K)^\times\to \Cl_0(\Delta)\to 1,
\end{equation}
where $(\OO_K^\times)_f$ is the image of the multiplicative group~$\OO_K^\times$ in ${(\OO_K/f\OO_K)}$. 

The group ${(\Z/f\Z)^\times(\OO_K^\times)_f}$  is ``not much bigger'' than   $(\Z/f\Z)^\times$. Precisely,
$$
[(\Z/f\Z)^\times(\OO_K^\times)_f:(\Z/f\Z)^\times]=[\OO_K^\times:\OO^\times]= 
\begin{cases}
2 &\text{if $D=-4$, $f>1$},\\
3&\text{if $D=-3$, $f>1$},\\
1&\text{otherwise}. 
\end{cases}
$$
An easy consequence is the following  formula for $h(\Delta)$:
\begin{equation}
\label{eclnum}
h(\Delta)= \frac{fh(D)}{[\OO_K^\times:\OO^\times]}\prod_{p\mid f}\left(1-\left(\frac Dp\right)p^{-1}\right),
\end{equation}
where ${(D/\cdot)}$ is the Kronecker symbol. 

\subsection{Orders with Class Groups Annihilated by~$2$}
\label{ssgen}

In this subsection we recall the famous result of Weinberger about imaginary quadratic orders whose class group is annihilated by~$2$. 
For a multiplicatively written abelian group~$G$ we denote by~$G^2$ its subgroup of squares: ${G^2=\{g^2: g\in G\}}$.

The group ${\Cl(\Delta)/\Cl(\Delta)^2}$ is usually called the \textsl{genus group} of~$\Delta$. It is known to be isomorphic to $(\Z/2\Z)^\mu$, where ${\mu=\mu(\Delta)\in \{\omega(\Delta)-1,\omega(\Delta)\}}$ and ${\omega(\cdot)}$ denote the number of distinct prime divisors. We may also remark that ${\mu(\Delta)= \omega(\Delta)-1}$ when ${\Delta=D}$ (and ${f=1}$). 

Already Euler studied discriminants~$\Delta$ with the property 
\begin{equation}
\label{ecl1}
|\Cl(\Delta)^2|=1,
\end{equation}
or, equivalently, ${\Cl(\Delta)\cong (\Z/2\Z)^\mu}$. (Of course, he used a different terminology.)
Chowla proved that the set of such~$\Delta$ is finite. Using a deep result of Tatuzawa~\cite{Ta51} about Siegel's zero, Weinberger~\cite{We73} improved on this, 
by showing that \textsl{field discriminants}~$D$ with this property are bounded explicitly with at most one exception. 

To state Weinberger's result precisely, denote by~$D'$ the square-free part of~$D$:
$$
D'=
\begin{cases}
D,& D\equiv 1 \bmod 4,\\
D/4, & D\equiv 0 \bmod 4.
\end{cases}
$$
\begin{proposition}[Weinberger \cite{We73}, Theorem~1] 
\label{pwein1}
There exists a negative integer~$D^\ast$ such that for any discriminant~$D$ of an imaginary quadratic field the property ${|\Cl(D)^2|=1}$ implies ${|D'|\le 5460}$ or ${D=D^\ast}$. 
\end{proposition}

Using existing methods~\cite{Wa04}, it is easy to determine the full list of~$D$ with ${|\Cl(D)^2|=1}$ and ${|D'|\le 5674}$. Since the group $\Cl(D)$ is a quotient of the group $\Cl(Df^2)$, if the latter has property~\eqref{ecl1} then the former does. Also, for every given~$D$ it is easy to find all possible~$f$ such that ${|\Cl(Df^2)^2|=1}$, using the description of the group $\Cl(Df^2)$ given in Subsection~\ref{ssclgr}. Hence the couples $(D,f)$ for which 
${|\Cl(Df^2)^2|=1}$ and ${|D'|\le 5674}$ can be  easily listed as well. 
This list is widely available in the literature since long ago; we reproduce it in Table~\ref{talist1}.
\begin{table}
\caption{Known~$\Delta$ with ${|\Cl(\Delta)^2|=1}$}
\label{talist1}
{\footnotesize
\begin{equation*}
\begin{array}{r|l}
h(\Delta)=1
&-3, -3\cdot2^2, -3\cdot3^2,-4, -4\cdot2^2,-7, -7\cdot2^2,-8,-11,-19,-43,-67,-163\\
\hline
h(\Delta)=2
&-3\cdot4^2, -3\cdot5^2, -3\cdot7^2,-4\cdot3^2, -4\cdot4^2, -4\cdot5^2,-7\cdot4^2,-8\cdot2^2, -8\cdot3^2,-11\cdot3^2,\\
&-15, -15\cdot2^2,-20,-24,-35,-40,
-51,
-52,-88,-91,
-115,-123,-148,\\
&-187,-232,-235,
-267,-403,-427\\
\hline
h(\Delta)\ge 4
&-3\cdot8^2,
 -7\cdot8^2,
  -8\cdot6^2,
  -15\cdot4^2, -15\cdot8^2,
 -20\cdot3^2,
 -24\cdot2^2,
 -35\cdot3^2,
 -40\cdot2^2,\\
&-84,
 -88\cdot2^2,
-120, -120\cdot2^2,
-132,
-168, -168\cdot2^2,
-195,-228,
 -232\cdot2^2,\\
&-280, -280\cdot2^2,
-312, -312\cdot2^2,
-340,
-372,
-408, -408\cdot2^2,
-420,
-435,\\
&-483,
-520, -520\cdot2^2,-532,
-555,
-595,
-627,
-660,
-708,
-715,\\
&-760, -760\cdot2^2,
-795,
-840, -840\cdot2^2,
-1012,
-1092,
-1155,
-1320, -1320\cdot2^2,\\
&-1380,
-1428,
-1435,
-1540,
-1848, -1848\cdot2^2,
-1995,
-3003,
-3315,
-5460
\end{array}
\end{equation*}  
}
\end{table}

It follows that Weinberger's result has the following consequence.

\begin{corollary}
\label{cwein1}
There exists a negative integer~$D^\ast$ such that  ${|\Cl(Df^2)^2|=1}$ implies that either ${\Delta=Df^2}$ appears in Table~\ref{talist1} or ${D=D^\ast}$.  
\end{corollary}

\begin{remark}
\label{rlist}
Class numbers of discriminants from Table~\ref{talist1}  are at most~$16$, and  the results of~\cite{Wa04}  imply that Table~\ref{talist1}  contains all~$\Delta$ with ${|\Cl(\Delta)^2|=1}$ and ${h(\Delta)\le 64}$. Hence if~$\Delta$ satisfies ${|\Cl(\Delta)^2|= 1}$ but does not appear in Table~\ref{talist1} then we must have ${h(\Delta)\ge 128}$.   

In particular, the first two lines of Table~\ref{talist1} give full lists of negative quadratic discriminants~$\Delta$ with ${h(\Delta)=1}$ and~$2$. 
\end{remark}

\section{Ring Class Fields and Complex Multiplication}
\label{sricf}

Let~$K$ be an imaginary quadratic field, and~$\OO$ an order in~$K$ of discriminant ${\Delta =Df^2}$. One associates to~$\OO$ an abelian extension  of~$K$ with Galois group $\Cl(\OO)$, called the \textsl{ring class field} of~$\OO$. We will denote it by $\ricf(\OO)$, or $\ricf(\Delta)$, or $\ricf (K,f)$. The canonical isomorphism ${\Cl(\OO)\to\gal(\ricf(\OO)/K)}$ is called the \textsl{Artin map}. For the details see, for instance, \cite[Section~9]{Co89}. 

The correspondence ${\OO\leftrightarrow\ricf(\OO)}$ is functorial in the following sense: if~$\OO'$ is a sub-order  of~$\OO$ then ${\ricf(\OO')\subset \ricf(\OO)}$, and we have the commutative diagram 
$$
\begin{array}{ccc}
\Cl(\OO)&\to&\gal(\ricf(\OO)/K)\\
\downarrow&&\downarrow\\
\Cl(\OO')&\to&\gal(\ricf(\OO')/K)
\end{array}
$$
where the horizontal arrows denote Artin maps and the vertical arrows are the natural maps of the class groups and the Galois groups. It follows that the Galois group of $\ricf(\OO)$ over $\ricf(\OO')$ is isomorphic to the kernel of  ${\Cl(\OO)\to \Cl(\OO')}$. In particular, $\gal(\ricf(\OO)/\ricf(\OO_K))$ is $\Cl_0(\Delta)$, the group introduced in Subsection~\ref{ssclgr}. (One may notice that ${\ricf(\OO_K)}$ is nothing else but the \textsl{Hilbert class field} of~$K$.)

\subsection{Compositum of Ring Class Fields}

In this subsection it will be more convenient to use the ``conductor notation'' $\ricf(K,f)$. 

As we have seen above, if ${f_1\mid f}$ then ${\ricf(K,f_1)\subset\ricf(K,f)}$. 
It follows that the compositum of two ring class fields $\ricf(K,f_1)$ and $\ricf(K,f_2)$ is a subfield of $\ricf(K,f)$, where ${f=\LCM(f_1,f_2)}$. It turns out that this compositum is ``almost always''  equal to $\ricf(K,f)$, but there are some exceptions.
Here is the precise statement. It is certainly known, but we did not find it in the available literature. 

\begin{proposition}
\label{plcm}
Let~$K$ be an imaginary quadratic field of discriminant~$D$ and $f_1,f_2$ positive integers. Set ${f=\LCM(f_1,f_2)}$. Then we have the following.

\begin{enumerate}

\item
If ${D\ne -3,-4}$ then ${\ricf(K,f_1)\ricf(K,f_2)=\ricf(K,f)}$. 

{\sloppy

\item
Assume that ${D\in \{-3-4\}}$. Then ${\ricf(K,f_1)\ricf(K,f_2)=\ricf(K,f)}$ either when one of $f_1,f_2$ is~$1$ or when  ${\gcd(f_1,f_2)>1}$. On the contrary, when ${f_1,f_2>1}$ and ${\gcd(f_1,f_2)=1}$,  the compositum ${\ricf(K,f_1)\ricf(K,f_2)}$ is a subfield of~$\ricf(K,f)$ of degree~$2$ for ${D=-4}$ and of degree~$3$ for ${D=-3}$.

}

\end{enumerate}

\end{proposition}

\paragraph{Proof (a sketch).}
To simplify the notation, we set 
\begin{align*}
L_0&=\ricf(K,1), & L_1&=\ricf(K,f_1), & L_2&=\ricf(K,f_2),\\
L&=\ricf(K,f),  & L'&=L_1L_2.  
\end{align*}
The mutual position of these fields is illustrated here:
$$
\begin{array}{ccccc}
&&L_1&&\\
&\nearrow&\downarrow&\searrow&\\
L_0&\longrightarrow&L'&\longrightarrow&L\\
&\searrow&\uparrow&\nearrow&\\
&&L_2&&
\end{array};
$$
We want to determine the degree ${[L:L']}$. 

We have ${\gal(L/L_0)=\Cl_0(Df^2)}$.  By~\eqref{ecl0structure}, this implies 
$$
\gal(L/L_0)= (\OO_K/f\OO_K)^\times/(\Z/f\Z)^\times(\OO_K^\times)_f,
$$
Similarly,
$$
\gal(L_i/L_0)= (\OO_K/f_i\OO_K)^\times/(\Z/f_i\Z)^\times(\OO_K^\times)_{f_i} \qquad (i=1,2).
$$
The Galois group ${\gal(L/L_i)}$ is the kernel of the natural map 
$$
(\OO_K/f\OO_K)^\times/(\Z/f\Z)^\times(\OO_K^\times)_f \stackrel{\pi_i}\to (\OO_K/f_i\OO_K)^\times/(\Z/f_i\Z)^\times(\OO_K^\times)_{f_i}.
$$
Hence $\gal(L/L')$ is the common kernel of the maps~$\pi_1$ and~$\pi_2$. It follows that ${\gal(L/L')=G/(\Z/f\Z)^\times(\OO_K^\times)_f}$, where~$G$ is the subgroup of $(\OO_K/f\OO_K)^\times$ consisting of ${x\in (\OO_K/f\OO_K)^\times}$ satisfying 
\begin{equation}
\label{egrg}
x\in (\Z/f\Z)(\OO_K^\times)_f \bmod f_i \qquad (i=1,2). 
\end{equation}
In particular, ${[L:L']=[G:(\Z/f\Z)^\times(\OO_K^\times)_f]}$.

If ${D\ne -3,-4}$ then  ${\OO_K^\times=\{\pm1\}}$, which implies that 
$$
G=(\Z/f\Z)^\times(\OO_K^\times)_f=(\Z/f\Z)^\times
$$ 
and ${L=L'}$. 

Now assume that ${D=-4}$. Then ${\OO_K^\times=\{\pm1,\pm \sqrt{-1}\}}$. When ${f_1=1}$ or ${f_2=1}$ the statement is trivial, so we may assume that ${f_1,f_2>1}$. Condition~\eqref{egrg} can be re-written as
{\small
\begin{equation}
\label{efour}
x\in \Z/f\Z\cup(\Z/f\Z)\sqrt{-1}\cup \bigl((\Z/f\Z)f_1+(\Z/f\Z)f_2\sqrt{-1}\bigr)\cup \bigl((\Z/f\Z)f_2+(\Z/f\Z)f_1\sqrt{-1}\bigr). 
\end{equation}}
If ${\gcd(f_1,f_2)>1}$ then the last two sets in~\eqref{efour} have no common elements  with $(\OO_K/f\OO_K)^\times$. We obtain
$$
G=\bigl(\Z/f\Z\cup(\Z/f\Z)\sqrt{-1}\bigr)\cap (\OO_K/f\OO_K)^\times= (\Z/f\Z)(\OO_K^\times)_f,
$$
and ${L=L'}$. 

If ${\gcd(f_1,f_2)=1}$ then each of the last two sets in~\eqref{efour} has  elements  belonging to $(\OO_K/f\OO_K)^\times$ but not to $(\Z/f\Z)(\OO_K^\times)_f$; for instance ${f_1+f_2\sqrt{-1}}$ and ${f_2+f_1\sqrt{-1}}$, respectively (here we use the assumption ${f_1,f_2>1}$). Hence
${[G:(\Z/f\Z)(\OO_K^\times)_f]>1}$. On the other hand, if~$x$ and~$y$ belong to the last two sets in~\eqref{efour}, then ${xy\in \Z/f\Z}$ if they belong to the same set, and ${xy\in (\Z/f\Z)\sqrt{-1}}$ if they belong to distinct sets. This shows that ${[G:(\Z/f\Z)(\OO_K^\times)_f]=2}$, and hence ${[L:L']=2}$. This completes the proof in the case ${D=-4}$.

The case ${D=-3}$ is treated similarly. We omit the details. \qed

\subsection{Complex Multiplication}
\label{sscm}
Ring class fields are closely related to the Complex Multiplication. Let ${\tau\in K}$ with ${\Im\tau>0}$ be such that   ${\OO=\End\langle\tau,1\rangle}$ (where ${\langle\tau,1\rangle}$ is the lattice generated by~$\tau$ and~$1$); one says that~$\OO$ is the \textsl{complex multiplication order} of the lattice ${\langle\tau,1\rangle}$. 

The ``Main Theorem of Complex Multiplication'' asserts that $j(\tau)$ is an algebraic integer generating over~$K$ the ring class field $\ricf(\OO)$. In particular, ${[K(j(\tau)):K]=h(\OO)}$. In fact, one has more:
\begin{equation}
\label{ehtau}
[K(j(\tau)):K]=[\Q(j(\tau)):\Q]=h(\OO). 
\end{equation}
The proofs can be found in many sources; see, for instance, \cite[Section~11]{Co89}. 

Since ${\ricf(\Delta)=K(j(\tau))}$ is a Galois extension of~$K$, it contains, by~\eqref{ehtau}, all the $\Q$-conjugates of $j(\tau)$. It follows that $K(j(\tau))$ is Galois over~$\Q$; in particular, the Galois group ${\gal(K/\Q)}$ acts on ${\gal(K(j(\tau))/K)}$. 

\begin{proposition}
\label{pdih}
The Galois group ${\gal(K/\Q)}$ acts on ${\gal(K(j(\tau))/K)}$ ``dihedrally'': if~$\iota$ is the non-trivial element of ${\gal(K/\Q)}$, then we have  ${\sigma^\iota=\sigma^{-1}}$ for any ${\sigma\in \gal(K(j(\tau))/K)}$. 
\end{proposition}

For a proof see, for instance \cite[Lemma~9.3]{Co89}.


\begin{corollary}
\label{c22}
The following properties are equivalent.

\begin{enumerate}
\item
\label{iab}
The field $K(j(\tau))$ is abelian over~$\Q$.

\item
\label{i2}
The Galois group ${\gal(K(j(\tau))/K)}$ is annihilated by~$2$ (that is, isomorphic to ${\Z/2\Z\times\cdots\times\Z/2\Z}$).

\item
\label{i22}
The Galois group ${\gal(K(j(\tau))/\Q)}$ is annihilated by~$2$. 

\item
\label{igal}
The field $\Q(j(\tau))$ is Galois over~$\Q$.

\item
\label{iab1}
The field $\Q(j(\tau))$ is abelian over~$\Q$.

\end{enumerate}

\end{corollary}

\paragraph{Proof.}
The implications \ref{iab}$\Rightarrow$\ref{i2}$\Rightarrow$\ref{i22} follow from Proposition~\ref{pdih}. The implication \ref{i22}$\Rightarrow$\ref{igal} is trivial. To see the implication \ref{igal}$\Rightarrow$\ref{iab1}, just observe that~\eqref{ehtau}  implies the isomorphism ${\gal(K(j(\tau))/K)\cong\gal(\Q(j(\tau))/\Q)}$. Finally, the implication \ref{iab1}$\Rightarrow$\ref{iab} is again trivial. \qed

\bigskip

\subsection{The Conjugates of $j(\tau)$}
\label{ssconj}
Let~$\tau$ and~$\OO$ be as in Subsection~\ref{sscm}, and let~$\Delta$ be the discriminant of~$\OO$.   As we already mentioned in the beginning of Subsection~\ref{sscm}, $j(\tau)$ is an algebraic integer of degree $h(\Delta)$. It is well-known that the $\Q$-conjugates of $j(\tau)$ can be described explicitly. Below we briefly recall this description.

Denote by ${T=T_\Delta}$ the set of triples of integers $(a,b,c)$ such that 
\begin{equation*}
\begin{gathered}
\gcd(a,b,c)=1, \quad \Delta=b^2-4ac,\\
\text{either\quad $-a < b \le a < c$\quad or\quad $0 \le b \le a = c$}
\end{gathered}
\end{equation*}
\begin{proposition}
\label{pallconj}
All $\Q$-conjugates of $j(\tau)$  are given by 
\begin{equation}
\label{ejabc}
j\left(\frac{-b+\sqrt{\Delta}}{2a}\right), \qquad (a,b,c)\in T_\Delta. 
\end{equation}
In particular, ${h(\Delta)=|T_\Delta|}$.
\end{proposition}
For a proof, see, for instance, \cite[Theorem~7.7]{Co89}.

The following observation  will be crucial: in the set $T_\Delta$  
there exists exactly one triple  ${(a,b,c)}$ with ${a=1}$. This triple can be given explicitly: it is 
$$
\left(1,r_4(\Delta),\frac{r_4(\Delta)-\Delta}4\right),
$$
where ${r_4(\Delta)\in\{0,1\}}$ is defined by ${\Delta\equiv r_4(\Delta)\mod 4}$. The corresponding number $j(\tau)$, where 
$$
\tau=\frac{-r_4(\Delta)+\sqrt\Delta}2,
$$ 
will be called \textsl{the dominant $j$-value} of discriminant~$\Delta$. It is important for us that it is much larger in absolute value than all its conjugates.

\begin{lemma}
\label{ldom}
Let $j(\tau)$ be the dominant $j$-value of discriminant~$\Delta$, with ${|\Delta|\ge 11}$, and let ${j(\tau')\ne j(\tau)}$ be  conjugate to $j(\tau)$ over~$\Q$. Then ${|j(\tau')|\le 0.1|j(\tau)|}$. 
\end{lemma}

\paragraph{Proof.} Recall the inequality 
${\bigl||j(z)|-|q_z^{-1}|\bigr|\le 2079}$, 
where ${q_z=e^{2\pi i z}}$, for~$z$ belonging to the standard fundamental domain of $\SL_2(\Z)$ on the Poincaré plane \cite[Lemma~1]{BMZ13}. We may assume that ${\tau=(-r_4(\Delta)+\sqrt\Delta)/2}$ and ${\tau'=(-b+\sqrt\Delta)/2a}$ with ${a\ge 2}$. Hence ${|q_\tau|= e^{\pi\sqrt{|\Delta|}}\ge e^{\pi\sqrt{11}}>33506}$ and ${|q_{\tau'}|\le |q_\tau|^{1/2}}$. We obtain 
$$
\frac{|j(\tau')|}{|j(\tau)|}\le \frac{|q_\tau|^{1/2}+2079}{|q_\tau|-2079}\le \frac{33506^{1/2}+2079}{33506-2079}<0.1,
$$
as wanted. \qed

\bigskip

The minimal polynomial of $j(\tau)$ over~$\Z$ is  called the  \textsl{Hilbert class polynomial}\footnote{or, sometimes, \textsl{ring class polynomial}, to indicate that its root generates not the Hilbert class field, but the more general ring class field} of discriminant~$\Delta$; it indeed depends only on~$\Delta$ because its roots are the numbers~\eqref{ejabc}. We will denote it $H_\Delta(x)$.

\section{Comparing two CM-fields}
\label{stwo}


 In this section we study the field equality ${\Q(j(\tau_1))=\Q(j(\tau_2))}$. We distinguish two cases: ${\Q(\tau_1)\ne \Q(\tau_2)}$, when we obtain the complete list of all possibilities, and ${\Q(\tau_1)= \Q(\tau_2)}$, where we will see that~$\Delta_1$ and~$\Delta_2$ are ``almost the same''. Here for ${i=1,2}$   we denote by~$\Delta_i$  the discriminant of the ``complex multiplication order'' ${\OO_i=\End\langle\tau_i,1\rangle}$, and write ${\Delta_i=D_i f_i^2}$ with the obvious meaning of~$D_i$ and~$f_i$.

\subsection{The case ${\Q(\tau_1)\ne \Q(\tau_2)}$}
\label{sstwodist}
In this subsection we investigate the case when ${\Q(j(\tau_1))=\Q(j(\tau_2))}$, but the fields $\Q(\tau_1)$ and $\Q(\tau_2)$ are distinct. It turns out that this is a very strong condition, which leads to a completely explicit characterization of all possible cases.

\begin{theorem}
\label{thedan}
Let~$\tau_1$ and~$\tau_2$ be quadratic numbers such that  ${\Q(j(\tau_1))=\Q(j(\tau_2))}$, but ${\Q(\tau_1)\ne \Q(\tau_2)}$. Then both~$\Delta_1$ and~$\Delta_2$ appear in Table~\ref{talist1}.
\end{theorem}

\paragraph{Proof.}
Denote by~$L$ the field ${\Q(j(\tau_1))=\Q(j(\tau_2))}$. If~$L$ is a Galois extension of~$\Q$, then  the group 
${\Cl(\Delta_1)=\Cl(\Delta_2)}$ is annihilated by~$2$ 
by Corollary~\ref{c22}, and we can use Corollary~\ref{cwein1}. Since ${D_1\ne D_2}$, at least one of the two discriminants~$D_1$ and~$D_2$ is distinct from~$D^\ast$; say, ${D_1\ne D^\ast}$. Then~$\Delta_1$ is in Table~\ref{talist1}. Since ${h(\Delta_1)=h(\Delta_2)}$, Remark~\ref{rlist} implies that~$\Delta_2$ is in Table~\ref{talist1} as well. (This argument goes back to Kühne \cite[Section~6]{Ku13}.)

Now assume that
\begin{equation}
\label{elngal}
\text{$L$ is not Galois over~$\Q$.}
\end{equation} 
We will show that this leads to a contradiction. Denote by~$M$ the Galois closure of~$L$ over~$\Q$; then   ${M=\Q(\tau_1,j(\tau_1))=\Q(\tau_2,j(\tau_2))}$. Define the Galois groups
$$
\begin{aligned}
&G=\gal(M/\Q), \qquad \tilN=\gal(M/\Q(\tau_1,\tau_2)),\\
&N_i=\gal(M/\Q(\tau_i))=\Cl(\Delta_i) \quad (i=1,2), 
\end{aligned}
$$
so that ${\tilN=N_1\cap N_2}$ and ${[N_1:\tilN]=[N_2:\tilN]=2}$. 

We claim the following: 
\begin{equation}
\label{eedan}
\text{the group~$\tilN$ is annihilated by~$2$}. 
\end{equation}
(One may mention that this is a special case of an observation  made independently by
Edixhoven~\cite{Ed98} and André~\cite{An98}.)

Indeed, let~$\iota_1$ be an element of~$G$ acting non-trivially on~$\tau_1$ but trivially on~$\tau_2$. 
Then for any ${\sigma\in N_1}$ we have ${\sigma^{\iota_1}=\sigma^{-1}}$ by Proposition~\ref{pdih}.   On the other hand, for any ${\sigma\in N_2}$ we have ${\sigma^{\iota_1}=\sigma}$. It follows that ${\sigma=\sigma^{-1}}$ for ${\sigma\in \tilN}$, proving~\eqref{eedan}. 

Thus, each of the groups~$N_1$ and~$N_2$ has a subgroup of index~$2$ isomorphic to $(\Z/2\Z)^\mu$ for some integer~$\mu$.  Hence each of~$N_1$ and~$N_2$ is isomorphic either to $(\Z/2\Z)^{\mu+1}$ or to ${\Z/4\Z\times(\Z/2\Z)^{\mu-1}}$. 

If, say, ${N_1\cong(\Z/2\Z)^{\mu+1}}$ then~$L$ is Galois over~$\Q$ by Corollary~\ref{c22}, contradicting~\eqref{elngal}. Therefore
$$
N_1\cong N_2\cong \Z/4\Z\times(\Z/2\Z)^{\mu-1}.
$$

Let ${\iota\in G}$ be the complex conjugation. Since~$i$ extends the non-trivial element of ${\Gal(\Q(\tau_1)/\Q)}$, the group ${H=\{1,\iota\}}$ acts on~$N_1$ dihedrally:  ${\sigma^\iota=\sigma^{-1}}$ for ${\sigma\in N_1}$, and we have  ${G=N_1\rtimes H}$. 
Let~$N_1'$ and~$N_1''$ be subgroups of~$N_1$ isomorphic to $\Z/4\Z$ and $(\Z/2\Z)^{\mu-1}$, respectively, such that ${N_1=N_1'\times N_1''}$.     Then~$H$ commutes with~$N_1''$ and acts dihedrally on~$N_1'$.   Hence
$$
G=N_1\rtimes H =(N_1'\rtimes H)\times N_1''\cong D_8\times (\Z/2\Z)^{\mu-1},
$$
where $D_{2n}$ denotes the dihedral group of $2n$ elements. 

Now observe that ${D_8\times (\Z/2\Z)^{\mu-1}}$ has only one subgroup isomorphic to the group ${\Z/4\Z\times(\Z/2\Z)^{\mu-1}}$; this follows, for instance, from the fact that both groups have exactly~$2^\mu$ elements of order~$4$. Hence ${N_1=N_2}$, which implies the equality ${\Q(\tau_1)=\Q(\tau_2)}$, a contradiction.\qed

\bigskip

Now one can go further and, inspecting all possible pairs of fields, produce the full list of number fields presented as  $\Q(j(\tau_1))$ and $\Q(j(\tau_2))$ with ${\Q(\tau_1)\ne \Q(\tau_2)}$. 

\begin{corollary}
\label{cojopa}
Let~$L$ be a number field with the following property: there exist  quadratic~$\tau_1$ and~$\tau_2$ such that ${L=\Q(j(\tau_1))=\Q(j(\tau_2))}$ but ${\Q(\tau_1)\ne \Q(\tau_2)}$. Then~$L$ is one of the fields in Table~\ref{tadouble}. 
\begin{table}
\caption{Fields presented as $\Q(j(\tau_1))$ and $\Q(j(\tau_2))$ with ${\Q(\tau_1)\ne \Q(\tau_2)}$}
\label{tadouble}
{\footnotesize
$$
\begin{array}{l|l|l|l}
\text{Field $L$}&[L:\Q]&\Delta&\Cl(\Delta)\\
\hline
\Q&1& -3, -4, -7, -8, -11, -12, -16, -19, -27, -28, -43,&\text{trivi-}\\ 
&&-67, -163&\text{al}\\
\hline
\Q(\sqrt{2})&2&-24, -32, -64, -88&\Z/2\Z\\
\Q(\sqrt{3})&2&-36, -48&\Z/2\Z\\
\Q(\sqrt{5})&2&-15, -20, -35, -40, -60, -75, -100, -115, -235&\Z/2\Z\\
\Q(\sqrt{13})&2&-52, -91, -403&\Z/2\Z\\
\Q(\sqrt{17})&2&-51, -187&\Z/2\Z\\
\hline
\Q(\sqrt{2},\sqrt{3})&4&-96, -192, -288&     (\Z/2\Z)^2\\
\Q(\sqrt{3},\sqrt{5})&4&-180, -240& (\Z/2\Z)^2\\
\Q(\sqrt{5},\sqrt{13})&4&-195, -520, -715& (\Z/2\Z)^2\\
\Q(\sqrt{2},\sqrt{5})&4&-120, -160, -280, -760&(\Z/2\Z)^2\\
\Q(\sqrt{5},\sqrt{17})&4&-340, -595&(\Z/2\Z)^2\\
\hline
\Q(\sqrt{2},\sqrt{3},\sqrt{5})&8&-480, -960&(\Z/2\Z)^3
\end{array}
$$
Explanations: 
\begin{enumerate}
\item
the third column contains the full list of discriminants~$\Delta$ of CM-orders $\End\langle\tau,1\rangle$ such that ${L=\Q(j(\tau))}$;
\item
the fourth column gives the structure of the class group $\Cl(\Delta)$ for any such~$\Delta$. 
\end{enumerate}
}\end{table}
\end{corollary}

\paragraph{Proof.}
This is just a calculation using \textsf{PARI}.\qed

\subsection{The case ${\Q(\tau_1)= \Q(\tau_2)}$}
\label{sstwoequ}

Now assume that ${\Q(j(\tau_1))=\Q(j(\tau_2))}$ and ${\Q(\tau_1)= \Q(\tau_2)}$. Denote by~$D$ be the discriminant of the number field ${\Q(\tau_1)= \Q(\tau_2)}$ and write 
${\Delta_i=f_i^2D}$ for ${i=1,2}$. 

\begin{proposition}
\label{pff2}
Assume that ${\Q(j(\tau_1))=\Q(j(\tau_2))}$ and ${\Q(\tau_1)= \Q(\tau_2)}$. Then either 
\begin{equation}
\label{eonetwohalf}
f_1/f_2\in \{1,2,1/2\} 
\end{equation}
or ${D=-3}$ and ${f_1,f_2\in\{1,2,3\}}$ (in which cases ${\Q(j(\tau_1))=\Q(j(\tau_2))=\Q}$). 
\end{proposition}

\paragraph{Proof.}
Put ${f=\LCM(f_1,f_2)}$. When ${D\ne -3,-4}$, Proposition~\ref{plcm} implies that 
\begin{equation}
\label{ecompos}
h(f^2D)= [\Q(\sqrt D,j(\tau_1),j(\tau_2)):\Q(\sqrt D)]
\end{equation}  
Since $j(\tau_1)$ and $j(\tau_2)$ generate the same field, we obtain 
\begin{equation}
\label{ehhh}
h(f_1^2D)=h(f_2^2D)=h(f^2D).
\end{equation}
Using~\eqref{eclnum} and~\eqref{ehhh}, we obtain
$$
\frac f{f_1}\prod_{\genfrac{}{}{0pt}{}{p\mid f}{p\nmid f_1}}\left(1-\left(\frac Dp\right)p^{-1}\right)=1,
$$
which implies that~${f/f_1\in \{1,2\}}$. Similarly, ${f/f_2\in \{1,2\}}$. Hence we have~\eqref{eonetwohalf}. 

Now assume that ${D\in \{-3,-4\}}$. If ${\gcd(f_1,f_2)>1}$ then we again have~\eqref{ecompos}, and the same argument proves~\eqref{eonetwohalf}.

If, say, ${f_1=1}$ then  either ${D=-4}$ and ${f_2\in \{1,2\}}$, in which case we again have~\eqref{eonetwohalf}, or 
${D=-3}$ and ${f_2\in \{1,2,3\}}$. 

Finally, assume that ${f_1,f_2>1}$ and ${\gcd(f_1,f_2)=1}$. Then ${f=f_1f_2}$ and Proposition~\ref{plcm} implies that 
\begin{equation}
\label{ehhlh}
h(f_1^2D)=h(f_2^2D)=\ell^{-1}h(f^2D),
\end{equation}
where ${\ell=2}$ for ${D=-4}$ and ${\ell=3}$ for ${D=-3}$. Using~\eqref{eclnum} and~\eqref{ehhlh}, we obtain
$$
f_i\prod_{p\mid f_i}\left(1-\left(\frac Dp\right)p^{-1}\right)=\ell \qquad (i=1,2),
$$
and a quick inspection shows that in this case ${D=-3}$ and ${\{f_1,f_2\}=\{2,3\}}$. \qed


\section{Proof of Theorem~\ref{thmain}}
\label{sproof}

We assume that ${P=(j(\tau_1),j(\tau_2))}$ belongs to a non-special straight line~$\ell$ defined over~$\Q$, and show that it satisfies either~\eqref{edeg1} or~\eqref{edeg2}. We define ${\Delta_i=D_if_i^2}$ as in the beginning of Section~\ref{stwo}.

Let ${A_1x_1+A_2x_2+B=0}$ be the equation of~$\ell$. Since~$\ell$ is not special, we have ${A_1A_2\ne 0}$, which implies, in particular, that 
\begin{equation}
\label{efields}
\Q(j(\tau_1))=\Q(j(\tau_2)).
\end{equation} 
We set
$$
L=\Q(j(\tau_1))=\Q(j(\tau_2)), \qquad h=h(\Delta_1)=h(\Delta_2)=[L:\Q].
$$
If   ${j(\tau_1),j(\tau_2)\in \Q}$ (that is, ${h=1}$) we are done.  From now on assume that 
$$
j(\tau_1),j(\tau_2)\notin \Q.
$$
If ${j(\tau_1)=j(\tau_2)}$ then~$\ell$  is the special line ${x_1=x_2}$, because it passes through the points $(j(\tau_1),j(\tau_1))$ and through all its conjugates over~$\Q$. Hence
\begin{equation}
\label{ene}
j(\tau_1)\ne j(\tau_2).
\end{equation}
If ${h=2}$ then we have~\eqref{edeg2}.  From now on assume that
\begin{equation}
\label{eh3}
h\ge 3, 
\end{equation}
and, in particular, 
\begin{equation}
\label{edelta16}
|\Delta_1|, |\Delta_2|\ge 23.
\end{equation}
We will show that this leads to a contradiction. 

\begin{remark}
\label{reasy}
Before proceeding with the proof, remark that, for a given pair of distinct discriminants~$\Delta_1$ and~$\Delta_2$ it is easy to verify whether there exists  a point ${(j(\tau_1),j(\tau_2))}$ on a non-special straight line defined over~$\Q$, such that~$\Delta_i$ is the discriminant of the CM-order ${\End\langle\tau_i,1\rangle}$. Call two polynomials ${f(x),g(x) \in \Q[x]}$ \textsl{similar} if there exist ${\alpha,\beta, \lambda\in \Q}$ with ${\alpha\lambda\ne 0}$ such that ${f(\alpha x+\beta)=\lambda g(x)}$. Now, a point ${(j(\tau_1),j(\tau_2))}$ as above exists if and only if the class polynomials $H_{\Delta_1}$ and $H_{\Delta_2}$ (see end of Subsection~\ref{ssconj}) are similar. This can be easily verified using, for instance, the \textsf{PARI} package. 
\end{remark}

\subsection{Both Coordinates are Dominant}

It turns out that we may assume, without a loss of generality, that both $j(\tau_1)$ and $j(\tau_2)$ are the \textsl{dominant  $j$-values} of corresponding discriminants, as defined in Subsection~\ref{ssconj}.

\begin{lemma}
\label{lbothdom}
Assume that~$\ell$ is a non-special  straight line containing a CM-point ${P=(j(\tau_1),j(\tau_2))}$ satisfying~\eqref{efields},~\eqref{ene} and~\eqref{eh3}. Then~$\ell$ contains a CM-point ${P'=(j(\tau_1'),j(\tau_2'))}$, conjugate to~$P$ over~$\Q$ and such that  $j(\tau_i')$ is the dominant $j$-value of discriminant~$\Delta_i$ for ${i=1,2}$. 
\end{lemma}

\paragraph{Proof.}
Since~$\ell$ is defined over~$\Q$, all  $\Q$-conjugates of~$P$ belong to~$\ell$ as well.   Replacing~$P$ by a $\Q$-conjugate point, we may assume  that $j(\tau_1)$ is the dominant $j$-value for the discriminant~$\Delta_1$. If $j(\tau_2)$ is the dominant value for~$\Delta_2$ we are done; so assume it is not, and show that this leads to a contradiction.

Since both $j(\tau_1)$ and $j(\tau_2)$ generate the same field of degree~$h$ over~$\Q$, the Galois orbit of~$P$ (over~$\Q$) has exactly~$h$ elements; moreover, each conjugate of $j(\tau_1)$ occurs exactly once as the first coordinate of a point in the orbit, and each conjugate of $j(\tau_2)$ occurs exactly once as the second coordinate.

It follows that there is a conjugate point $P^\sigma$  such that the second coordinate $j(\tau_2)^\sigma$  is the  dominant $j$-value for~$\Delta_2$; then its first coordinate $j(\tau_1)^\sigma$ is  not dominant for~$\Delta_1$ because ${P^\sigma\ne P}$. Since ${h\ge 3}$, there exists yet another point~$P^{\sigma'}$ with both coordinates  not dominant for the respective discriminants.

All three points~$P$,~$P^\sigma$ and~$P^{\sigma'}$ belong to~$\ell$. Hence
$$
\begin{vmatrix}
1& j(\tau_1)&j(\tau_2)\\
1& j(\tau_1)^\sigma&j(\tau_2)^\sigma\\
1& j(\tau_1)^{\sigma'}&j(\tau_2)^{\sigma'}
\end{vmatrix}=0.
$$
The determinant above is a sum of~$6$ terms: the ``dominant term'' ${j(\tau_1)j(\tau_2)^\sigma}$ and~$5$ other terms. Each of the other terms is at most ${0.1|j(\tau_1)j(\tau_2)^\sigma|}$  in absolute value: this follows from Lemma~\ref{ldom}, which applies here due to~\eqref{edelta16}.  Hence the determinant cannot vanish, a contradiction. \qed

\subsection{Completing the proof}
After this preparation, we are ready to complete the proof of Theorem~\ref{thmain}. Thus, let 
${P=(j(\tau_1),j(\tau_2))}$ belong to a straight line~$\ell$ defined over~$\Q$. We assume that~\eqref{efields},~\eqref{ene} and~\eqref{eh3} are satisfied, and we may further assume that $j(\tau_1)$, $j(\tau_2)$  are the dominant $j$-values of~$\Delta_1$,~$\Delta_2$, respectively.

If ${\Q(\tau_1)\ne \Q(\tau_2)}$ then Corollary~\ref{cojopa} applies, and all possible~$L$,~$\Delta_1$ and~$\Delta_2$ can be found in Table~\ref{tadouble}. In particular,  we have only 6 possible fields~$L$ and 15 possible couples $\Delta_1,\Delta_2$. All of the latter are ruled out by verifying (using \textsf{PARI}) that the corresponding Hilbert class polynomials $H_{\Delta_1}(x)$ and $H_{\Delta_2}(x)$ are not similar, as explained in Remark~\ref{reasy}.

If ${\Q(\tau_1)= \Q(\tau_2)}$ then Proposition~\ref{pff2} applies, and we have ${f_1/f_2\in\{1,2,1/2\}}$. 
Since both $j(\tau_1)$ and $j(\tau_2)$ are dominant, the case ${f_1=f_2}$ is impossible: there is only one dominant $j$-value for every given discriminant, and we have ${j(\tau_1)\ne j(\tau_2)}$ by~\eqref{ene}. Thus, ${f_1/f_2\in \{2, 1/2\}}$.

Assume, for instance, that ${f_2=2f_1}$. Write ${\Delta_2=\Delta}$, so that ${\Delta_1=4\Delta}$. Since both $j(\tau_1)$ and $j(\tau_2)$ are dominant, we may choose
$$
\tau_1= \frac{-r_4(4\Delta)+\sqrt{4\Delta}}{2}=\sqrt\Delta, \qquad \tau_2= \frac{-r_4(\Delta)+\sqrt{\Delta}}{2}
$$
It follows that ${\tau_2=\frac12\gamma(\tau_1)}$, where 
$$
\gamma=\begin{pmatrix}
1&-r_4(\Delta)\\0&1
\end{pmatrix}\in \SL_2(\Z). 
$$
Hence the point ${P=(j(\tau_1),j(\tau_2))}$ belongs to the modular curve $Y_0(2)$ realized as the plane curve ${\Phi_2(x_1,x_2)=0}$, where 
\begin{align*}
\Phi_2(x_1,x_2) =& -x_1^2x_2^2 + x_1^3 + x_2^3 + 1488 x_1^2x_2 +1488 x_1x_2^2 +40773375x_1x_2 \\
&- 162000x_1^2-162000x_2^2+ 8748000000x_1+8748000000x_2\\
&-157464000000000
\end{align*}
is the modular polynomial of level~$2$. 
Since ${\deg \Phi_2=4}$ and~$P$ belongs to a straight line over~$\Q$, the coordinates of~$P$ generate a field of degree at most~$4$ over~$\Q$. Thus, ${3\le h\le4}$.

Looking into existing class number tables (or using \textsf{PARI}) one finds that there exist only 5 negative discriminants~$\Delta$ such that 
${h(\Delta)=h(4\Delta)\in \{3,4\}}$:
\begin{align*}
h=3:&\quad -23, -31; \\
h=4:&\quad -7\cdot3^2, -39,-55.
\end{align*}
Verifying that the  polynomials $H_\Delta$ and $H_{4\Delta}$ for these values of~$\Delta$ are not similar is an easy calculation with \textsf{PARI}. \qed

\begin{remark}  We conclude the article with some computational remarks. 
\label{rall}
\begin{enumerate}
\item
As indicated in the introduction, it is very easy to list all CM-points satisfying~\eqref{edeg1} or~\eqref{edeg2}; call them \textsl{rational} and \textsl{quadratic}, respectively.  

There exist exactly~$13$ discriminants~$\Delta$ with ${h(\Delta)=1}$,  the first line  of Table~\ref{talist1} lists them all. Hence there exist~$169$ rational CM-points, and listing them explicitly is plainly straightforward. 

As for the quadratic CM-points, there are two kinds of them:
points  with 
\begin{equation}
\label{econj}
\Delta_1=\Delta_2=\Delta, \quad h(\Delta)=2, \quad \text{$j(\tau_1)$ and $j(\tau_2)$ are conjugate over~$\Q$},
\end{equation}
and 
points  with 
\begin{equation}
\label{enonconj}
\Delta_1\ne\Delta_2, \quad h(\Delta_1)=h(\Delta_2)=2, \quad \Q(j(\tau_1))=\Q(j(\tau_2)). 
\end{equation}
There exist exactly~29 negative quadratic discriminants~$\Delta$ with ${h(\Delta)=2}$, see the second line of Table~\ref{talist1} for the complete list. Hence there exist~$29$, up to conjugacy, quadratic CM-points satisfying ~\eqref{econj}. 

The quadratic CM-points satisfying~\eqref{enonconj} can be extracted from the ``quadratic'' part of  Table~\ref{tadouble}. Indeed, if ${\Q(\tau_1)\ne\Q(\tau_2)}$ then $\Delta_1,\Delta_2$ are in Table~\ref{tadouble} by Corollary~\ref{cojopa}. And if ${\Q(\tau_1)=\Q(\tau_2)}$ then ${\Delta_1/\Delta_2\in \{4,1/4\}}$ by Proposition~\ref{pff2}. Inspecting the list of the~29 discriminants with ${h(\Delta)=2}$, we find that the only possibility is ${\{\Delta_1,\Delta_2\}=\{-15,-60\}}$. Both these values appear in Table~\ref{tadouble} in the line corresponding to the field $\Q(\sqrt5)$. 

Looking into Table~\ref{tadouble} we  find that there exist 
$$
4(4-1)+2(2-1)+9(9-1)+3(3-1)+2(2-1)=94
$$
(ordered) pairs $(\Delta_1,\Delta_2)$ as in~\eqref{enonconj}. Each pair gives rise to~two, up to conjugacy, points satisfying~\eqref{enonconj}. So, up to conjugacy, there are 188 points satisfying~\eqref{enonconj}, and ${188+29=217}$ quadratic CM-points altogether. Again, listing them explicitly is a straightforward computation. 

\item
Thomas Scanlon (private communication) asked whether there exists a non-special straight line over~$\C$ passing through more than~$2$ CM-points. Since 
$$
\det\begin{bmatrix}
1728& -884736000\\
287496& -147197952000 
\end{bmatrix}=0,
$$
the points $(0,0)$, $(1728,287496)$ and $(-884736000,-147197952000)$ belong to the same straight line, and so do the points $(0,0)$, $(1728,-884736000)$ and $(287496,-147197952000)$. 
Notice that 
\begin{align*}
&j\left(\frac{-1+\sqrt{-3}}{2}\right)=0, \quad j(\sqrt{-1})=1728, \quad j(2\sqrt{-1})=287496,\\
&j\left(\frac{-1+\sqrt{-43}}{2}\right)=-884736000, 
\quad j\left(\frac{-1+\sqrt{-67}}{2}\right)=-147197952000. 
\end{align*}
We verified that  (up to switching the variables $x_1,x_2$) these are the only such lines defined over~$\Q$. Precisely:

\begin{itemize}
\item
no~$3$ rational CM-points, with the exceptions indicated above, lie on the same non-special line;

\item
no line passing through conjugate quadratic CM-points contains a rational CM-point;

\item
lines defined by pairs of conjugate quadratic CM-points are all pairwise distinct.

\end{itemize}
The verification is a quick calculation with \textsf{PARI}.

It is not clear to us whether the examples above (and absence of other examples) are just accidental, or admit some conceptual explanations. 

\item
All the computations for this paper take about 20 seconds. 
\end{enumerate}
\end{remark}

{\footnotesize

\begin{tabular}{lll}
\textbf{Bill Allombert} & \textbf{Yuri Bilu} & \textbf{Amalia Pizarro-Madariaga}\\
Université de Bordeaux &  Institut de Mathématiques de Bordeaux & Instituto de Matemáticas\\ 
IMB, UMR 5251 & Université de Bordeaux & Universidad de Valpara\'iso\\
F-33400 Talence, France & 351 cours de la Libération &Chile\\
&33405 Talence, France&\\
CNRS&&\\ IMB, UMR 5251&&\\
F-33400 Talence, France&&\\
&&\\
INRIA&&\\ F-33400 Talence, France
\end{tabular}

}

\end{document}